\documentstyle[12pt]{article}
\begin{document}
\def\A{{\cal A}}
\def\B{{\cal B}}
\def\C{{\cal C}}
\def\D{{\cal D}}
\def\E{{\cal E}}
\def\F{{\cal F}}
\def\G{{\cal G}}
\def\I{{\cal I}}
\def\J{{\cal J}}
\def\K{{\cal K}}
\def\M{{\cal M}}
\def\N{{\cal N}}
\def\Q{{\cal Q}}
\def\R{{\cal R}}
\def\T{{\cal T}}
\def\U{{\cal U}}
\def\V{{\cal V}}
\def\W{{\cal W}}
\def\X{{\cal X}}
\def\Y{{\cal Y}}
\def\Z{{\cal Z}}
\def\be{\begin{equation}}
\def\ee{\end{equation}}
\def\la{\label}
\def\sc{\scriptstyle}
\def\tf{\tilde\F}
\def\tfk{\tf^k}
\def\da{\dagger}
\def\al{\aleph}
\def\alk{\al_k}

\begin{center}
{\bf V.A.Buslov}\end{center}

\vskip1cm

\begin{center}
{\bf\Large  Hierarchy Structure of Graphs and Weighted Condensations
 }
\end{center}
\vskip 8mm

\centerline{Abstract}

\centerline{Buslov V.A.}

\parbox[t]{14cm}
{\small By natural way the hierarchy structure is introduced on directed
graphs with weighted adjacencies. Embedded system of algebras of subsets
of the set of vertices of such digraph and it's consolidations, which
vertices are the elementary sets of corresponding algebra, are
constructed. Weights of arcs of consolidated graphs are determined.}

\section{Introduction words}

The very first essential theorem on graphs (except Euler's (1736)
solution of K\"enigsberg bridges problem) was formulated by Kirchhoff
\cite{K} (1847) who considered graphs as networks of conducting wires. In
this theorem Kirchhoff computed the number of connected subgraphs
containing all vertices without circuits (spanning trees) for the aims of
analyzing electric chains.

The next step in investigating tree-like structure of graph was the
following. One can attach to every edge (unordered pair of vertices) some
quantity (weight) and ask a question how to find the spanning tree of such
graph (graph with weighted adjacencies) with minimal weight (having the
minimal sum of edge weights). To clarify the idea of weights and the
problem formulated one can use Kirchhoff's approach. Let us assume that
some number of points (vertices of graph) are connected by wires (edges
of graph) in an arbitrary manner. The resistance (or length) of wire
connecting two points is the weight of corresponding edge. The question is
which wires to remove and which to reserve in order to get such network
where all points are still connected by wires (may be through other
points) but the summarized resistance (or length) of remaining wires be
minimal.  Of course such network must be a spanning tree. This significant
problem is in graph theory one of few ones having a number of effective
algorithms.  According to one of them at the first step one takes an edge
of minimal weight as a first intermediate graph (if it is not a single
minimal weight edge one takes any of them), then the procedure is
recurrent one.  If $k$-th intermediate graph is constructed one have to
add to it the edge of minimal weight among all others such as resulting
graph does not contain circuits. At the $(N-1)$-th step, where $N$ is a
number of all vertices, one gets spanning tree of minimal weight.  So at
every step one gets a forest (graph without circuits), which connected
components are trees, and every step means that two of them are joined
into single tree by adding some edge connecting them. But for directed
graph (digraph), where edges replaced by arcs (ordered pairs of vertices),
the equivalent procedure was absent to our knowledge, and Theorem 1 of
present paper shows us the way of such constructing.

To make clear the idea of weights and the problem to solve in the
situation of digraphs we need to replace the classical picture of graph as
an electric network by another one. Let us consider the picture of
potential relief on a smooth manifold $M$, defined by some real smooth
function $V(x)$, $x\in M$.  The points of local minimum of the potential
we associate as vertices of digraph and the potential bar necessary to
overtake in order to leave point $x_i$ of local minimum and leave the
corresponding fundamental region $\Omega_i$ (the point $x_i$ attraction
region of dynamic system $\dot x=-\nabla V(x)$) and fall into another
fundamental region $\Omega_j$, which has a common boundary with
$\Omega_i$, we assume to be weight $V_{ij}$ of arc $(i,j)$ of digraph.

The main idea of constructing  spanning trees of minimal weight is
similar to the method above.  At $k$-th step we construct forest of
minimal weight among forests having $k$ arcs by connecting by arc two
trees of the previous step into a single tree but with serious addition.
One of these two trees must be reconstructed into new one before, and only
then we connect them.

Let us explain partially here the reason of this reconstructing. There are
two different orientation types of directed trees in the situation of
digraphs.  In one of them moving from any vertex along arcs according to
orientation one falls into a single vertex called a root (just this type
of directed trees we consider at this parer).  By changing the direction
of all arcs we get another type of directed trees. The complexity of
directed case can be demonstrated by the following fact: it is easy to
construct directed graph with weighted adjacencies such that any of
its spanning trees (independently of the orientation type) of minimal
weight does not contain the minimal weight arc!

Returning to our example of potential field on manifold the problem of
constructing spanning tree having minimal weight can be reformulated as
the next one. The problem is to chose some point of local minima (the
root of the tree -- it would be the point of global minimum) and to bind
points of local minima by one-side ways so that moving from any point
along these ways one gets into a the chosen point and the sum of
potential bars corresponding to these ways must be minimal.

But the main aim of present paper is not to construct directed spanning
trees (it is only a collateral result), although tree remains the
elementary brick of our construction. In fact we consider a sequence of
ordered grainings of given graph nested into each other. Every of these
grainings except the last one is represented by simpler graph (the
"descendant") still possessing essential properties of its "ancestor".

To illustrate the idea of grainings let us return to our example of
manifold. One can imagine the particle, moving on this manifold. If the
particle is at the moment at some point of local minima and possesses some
fixed additional energy and this energy is quite enough to overtake some
potential bars, the particle does not notice these bars and does not
distinguish the fundamental regions they (bars) connect.  So our manifold
divides into more rough regions than fundamental ones. Increasing the
level of this additional energy one gets more rough subdivision etc.

So coming back to graphs it is naturally to ask such a question.  Can we
look on graphs not in detail, but rather roughly without noticing
unessential connections and uniting the vertices themselves into some
grainings, among which to establish new connections?  Continuing this
logic we could look another time on such graining graph and to enlarge it
one more time and so on.  We suggest such kind of hierarchy approach here.
Inside the situation when arcs (for directed) or edges (for non-directed)
graphs are of no weights, the proposition about graph hierarchy structure
is rather poor. For usual (non-directed) graphs we can speak in such sense
only about connected components, for directed -- about graph of
condensations, which elements are the sets of mutually accessible
vertices. For directed (and even for non-directed) graphs with weighted
adjacencies the situation is far richer.  One could try to introduce by
force some decompositions on clusters of vertices and then has a problem
how to determine adjacencies connecting these clusters, but our approach
is natural one.  It implies that graph itself contains the whole
information about the number of hierarchically nested decompositions and
about clusters inside every decomposition and this is determined mainly
not by the number of vertices and not by the number of arcs (edges) but
mostly by forces of these adjacencies (by weights).  In this connection it
comes to light that one can neglect  the part of adjacencies without any
waste and such a neglect does not affect on clusters inside decomposition
and the number of decompositions either.

On the structure of work. At section 2 we give the necessary definitions
and designations.  At section 3 we introduce some new definitions required
to formulate the general method using in proofs. Next section is devoted
to investigating the properties of directed forests, which are the
factors of the initial graph and have the minimal weight among all forests
of $k$ trees under different $k$. It turns out that such extreme forests
allow to construct at section 5 nested system of algebras of subsets of
the set of vertices (they determine the hierarchy) and to investigate
their properties. At section 6 using results of the previous section we
construct some kind of enlarged graphs, which we call by weighted
condensations.

\section{Main definitions }

In graph theory the unification of designations and even terminology
proper have not complete yet. So let us give first necessary definitions
and designations.

Let $G$ be graph (non-directed). By $\V G$ and $\E G$ we denote the set of
its vertices and edges (unordered pairs of vertices) respectively.

Let $\X$ be non-empty set and $\X ^2$ -- its Cartesian square and let
$\U\subseteq\X^2$.  Pair $G=(\X ,\U)$ is called {\it directed graph
(digraph)}.  The elements of the set $\X$ are called {\it vertices} and
elements of the set $\U$ are called {\it arcs}. We use $\V G$ and $\A G$
to denote the set of vertices and arcs of $G$ respectively.

Let $a=(i,j)$ be an arc, vertices $i$ and $j$ are called  {\it an origin}
and {\it a terminus} of $a$ respectively.  The arc $(i,i)$ with
coinciding origin and terminus is called {\it a loop}. The number of arcs
coming out of (into) the vertex $i$ is called {\it outdegree}
$d^+(i)$ ({\it indegree} $d^-(i)$) of the vertex $i$ .

Digraph having several arcs with common origin and common terminus is
called {\it multidigraph} and such arcs are called {\it parallel}. If
every edge (arc) of (di)graph possesses some value (weight), such (di)graph
is called {\it (di)graph with weighted adjacencies}.

We use sometimes the term "graph" in wide sense designating by it digraphs
and multidigraphs with weighted adjacencies also if it is not lead to
misunderstanding.

Graph $H$ is called {\it subgraph} of graph $G$ if $\V H\subseteq\V G$, $
\A H\subseteq\A G$.  Subgraph $H$ is called {\it spanning} subgraph (or
{\it factor}) if $\V H=\V G$.  Subgraph $H$ is called {\it induced} (or
more completely -- subgraph induced by the set $\U\subset\V G$) if $\V H
=\U$ and $(i,j)\in\A H$ means that  $(i,j)\in\A G$ and $\{i,j\}\subset
\U$. We designate subgraph of $G$ induced by the set  $\U$ as $G|_\U$.

{\it Directed circuit of length M} is digraph with set of vertices
$\{ x_1$,$x_2$,$\cdots $,$ x_M\}$ and with arcs $(x_j,x_{j+1}) ,
j=1,2,\cdots M-1$ and $(x_M,x_1)$.

{\it Walk (noncyclical) of length $M-1$} is digraph with set of vertices
$\{ i_1$,$2$, $\cdots $,$M\} $ and with arcs $(i_j,i_{j+1})$, $i=1$,
$\cdots $,$M-1$ .  Such walk we designate $i_1\cdot i_M$-walk. {\it
Semiwalk of length $M-1$} is digraph with set of vertices $\{
i_1,i_2,\cdots ,M\}$ and its arcs are either $(i_j,i_{j+1})$ or
$(i_{j+1},i_j)$ , $i=1,2,\cdots ,M-1$. The vertex $j$ is said to be {\it
accessible} ({\it attainable}) from the vertex $i$ in graph $G$ if there
is $i\cdot j$-walk in $G$. Digraph is called {\it strong} (or {\it strong
connected}) if all its vertices are mutually attainable. Digraph $G$ is
called {\it weak} if for every pair of vertices there is a semiwalk
connecting them in $G$.

Any maximal with respect to including weak subgraph of graph $G$ is called
its connected component (or simply -- component). {\it Strong component}
of $G$ is any maximal with respect to including its strong subgraph.

Other  definitions we will cite as it is necessary.

\section{Other definitions, designations and pre\-de\-ter\-mi\-ned
operations}

Graph (non-directed) possessing no cycles is a {\it forest}. Connected
component of a forest is a {\it tree}. For trees with $\S $ as a set
of vertices we use the notation $T(\S  )$.

There are two kinds of directed forests at the situation of digraphs. Here
we call by {\it forest} the digraph without circuits, in which outdegree
of every vertex is equal to zero or to one ($d^+(i)=0,1$). Arcwise
connected components of forest are called {\it trees}. The only vertex $i$
of tree which outdegree is equal to zero ($d^+(i)=0$) is called {\it root}
of a tree. The set of roots of the forest $F$ we designate by $\W ^F
$.  The tree of $F$ with root $i$ we designate $T^F_i$.

Let $V$ be graph (directed or non-directed). We use notation $\F^k(V)$ for
the set of spanning forests having $k$ trees and being subgraph of $V$.

We call the vertex $i$ {\it rear} to the vertex $j$ in graph $G$, and
correspondingly $j$ {\it front} to $i$ if there is $i\cdot j$-walk in $G$.
{\it Front} ({\it rear}) {\it enclosing} of the vertex $i$ in graph $G$ is
the set of terminuses of arcs outcoming from (origins of arcs coming
into) $i$. For such set we use notation $\N ^+_G(i)$ ($\N ^-_G(i)$).

{\bf Remark.} If digraph is a forest, the front-rear relation is a
relation of partial ordering. Such definition implicates that vertices can
be connected by the rear-front relation  only if they are in the same tree
of the forest.  The root of the forest is the front vertex to all vertices
of the tree.

We say that in graph $G$ an arc comes out of the set $\U$ and comes into
the set $\V$, if there is at list one arc which origin belongs to $\U$
and the terminus belongs to $\V$ in graph $G$.  We also say that arc
comes out of the set $\U$, if there is at least one arc which origin
belongs to $\U$ but terminus does not.

If $\D$ is some subset of the set of all vertices of digraph $G$ we call
the set of terminuses (origins) of arcs outcoming from (coming into) the
set $\D$ as front (rear) enclosing of $\D$, which is naturally to
designate as $\N^+_G(\D)$ ($\N ^-_G(\D)$).

In the following we will prove the existence of forests having some special
properties. The general method of such proofs consists of sequent steps.
It is necessary to take two concrete graphs with the same set of vertices
and to select some subset $\D$ of vertices. Next, we exchange between each
other arcs outcoming from the vertices of $\D$ in this graphs and then we
investigate properties of the new graphs resulting in such exchange.
Thereby is naturally to introduce the following definition.

Let $F$ and $G$ be two graphs with the same set of vertices and the set
$\D$ is some subset of the set of vertices. We will say that the graph
$H$ is $\D$-{\it exchange} of  $F$ by $G$, if  $H$ is a result of
exchanging in graph $F$ arcs, outcoming from the vertices of the set $\D$,
onto arcs that outcome from these vertices in graph $G$.  Our interest
relates to the situation, where $F$ and $G$ are forests and in addition
$\D$-exchange of $F$ by $G$ and moreover at the same time $\D$-exchange of
$G$ by $F$ are forests too. So at first let us formulate the criterion of
$\D$-exchange to be a forest.

{\bf Criterion.} {\it Let $F$ and $G$ be two an arbitrary forests with
the set of vertices $\N =\{ 1,2,\cdots ,N\}$,  $\D $ -- some subset of $\N
$. Let $H$ is  $\D$-exchange $F$ by $G$.  Then  graph $H$ is a forest then
and only then, if every vertex $i\in\N^+_G\D$ is not rear in $F$ with
respect to those vertices from $\D$ which are rear to $i$ in $G$.}

{\bf Proof.} As any way starting from the vertex $i\in\N^+_G(\D)$ in graph
$F$ (and as a sequence in  $H$ too) by the condition can not include those
vertices of the set $\D$, which are rear to $i$ in $G$ (and as a sequence
in $H$), so $H$ does not contain circuits.  Further, not more than one arc
comes out from any vertex in $H$, so $H$ is a forest.

{\bf Sequence 1.} {\it Let $F$ -- be a forest and $\D$ be some subset of
the set of vertices, i) if there are not any arcs coming into $\D$ in
$F$, so $\D$-exchange of $F$ by any forest $G$ is a forest; ii) if there
are not any arcs coming out of $\D$ in $F$, so $\D$-exchange of any
forest $G$ by $F$ is a forest.}

{\bf Sequence 2.} {\it Let  $T^F$ be a tree of the forest $F$ and $T^G$
-- tree a of the forest $G$, and let $\D =\V T^F$ ($=\V T^F\cap\V T^G$,
$=\V T^F\setminus\V T^G$) then $\D$-exchange of $F$ by $G$ and $
\D$-exchange of $G$ by $F$ are forests.}

{\bf Sequence 3.} {\it Let $T^F$ be a tree of the forest $F$ and $T^G$
be a tree of the forest $G$, $\C =\V T^F\cap\V T^G$, and $\D\subseteq\C$,
such, that there are not arcs coming into $\D$ in $F$ and terminuses of
arcs coming from $\D$ do not belong to $\C$ , then $\D$-exchange of $F$
by $G$ and $\D$-exchange of $G$ by $F$ are forests.}

{\bf Sequence 4.} {\it  Let $T^F$ be a tree of the forest  $F$ and $T^G$
be a tree of $G$, $\C =\V T^G\setminus \V T^F$, and $\D$ -- be the set of
all vertices from $\C$, such as the walk starting from any of them in the
forest $G$ passes through the set $\V T^F$, then $\D$-exchange of $F$ by
$G$ and $\D$-exchange of $G$ by $F$ are forests.}

\section{Related forests}

Let $V$ be digraph with real weighted adjacencies $v_{ij}$ on the set of
vertices $\N =\{ 1,2,\cdots ,N\}$. We will consider factors of $F$ being
forests and containing of $k=1,2,\ldots ,N$ trees (the set of such forests
we designate $\F ^k(V)$).  Under {\it weight} $\Sigma^F$ of $F$ we
understand the following quantity:

$$
\Sigma^F=\sum_{(i,j)\in \A F}v_{ij} \ .
$$
The minimum of weight over all forests $F\in\F^k(V) $ consisting
exactly of $k$ trees we designate as $\varphi_k$:

$$ \varphi_k=\min_{F\in\F^k(V)}\Sigma^F \ .  $$
If $\F^k(V)=\emptyset$ we suppose  $\varphi_k=\infty$. In the following we
write $\F^k$ instead of $\F^k(V)$ in cases when it is clear subgraphs of
which graph $V$ are under consideration.

Let us pick out the subset $\tfk $ from the set of forests $\F^k$,
consisting of forests with the minimum weight:  $F\in\tfk \Leftrightarrow
\F\in\F^k$ and $ \ \Sigma^F=\varphi_k $. Such forests we call extreme.

Let us study extreme (giving minimum) forests from $\tfk (V)$ under
different $k$. It turns to be that they have some kind of "genetic" link.
In particular it is valid the following

{\bf Proposition 1.} {\it Let under some $k=1,2,\cdots ,N-1$, the set
$\F^k$ is not empty, then for any forest $F\in\tf^{k+1}$ there is at
least one $G\in\tfk$ (and for any forest $G\in\tfk$ there is  $F\in\tf^
{k+1}$) such that the set of vertices of any tree of the forest $F$ is
contained in the set of vertices of some tree of the forest $G$}.

{\bf Remark.} Just the formulation of this proposition means that as the
forest $G\in\tfk$ contains one tree less than "relative" to it forest
$F\in\tf^{k+1}$, so  the sets of vertices of $k-1$ trees of the forest $G$
coincide with the sets of vertices of corresponding trees of the forest
$F$, the set of vertices of the last tree of the forest $G$ is conjunction
of sets of vertices of last two trees of the forest $F$.

In actual we will prove more powerful fact. Preliminary we give one
definition.

Let us agree upon to call the forest $F\in\F^{k+1}$ with roots (exactness
to the numeration) $ \ 1,2,\cdots ,k+1 \ $ as an  {\it ancestor} of the
forest $G\in\F^k$ with roots $1,2,\cdots ,k$, and correspondingly the
forest $G\in\F^k$ to call as a {\it descendant} of the forest $F\in\F
^{k+1}$ if $T^F_i=T^G_i, \ \ i=1,2,\cdots ,k-1$, $T^F_k\subset T^G_k$, and
subgraph $G|_{\V T^F_{k+1}}$ of the forest $G$ (or, which is the same,
subgraph of the tree $T^G_k$) induced by the set $\V T^F_{k+1}$ is a tree
(under this it may coincide with the tree $T^F_{k+1}$ or not).

The following theorem  tells us on the minimum changes one must to provide
to get a forest belonging to the set $\tfk$ from a forest belonging to the
set $\tf^{k+1}$ and vice versa.

{\bf Theorem 1 (on "relatives").} {\it Let under some  $k=1,2,\cdots ,N-1$
the set $\F^k$ is not empty, then any forest $F\in\tf^{k+1}$ has a
descendant in the set $\tfk$ and any forest $G\in\tfk$ has an ancestor
in the set $\tf^{k+1}$.}

{\bf Proof.} Let us prove that any forest $F\in\tf^{k+1}$ has a descendant
in the set $\tfk$.  Let $F$ and $H$ be arbitrary forests from the sets
$\tf^{k+1}$ and $\tfk$ respectively. As the power of the set of roots
$\W^F$ of the forest $F$ is one unit more than the power $|\W^H|=k$, and
as in any forest not more than one arc goes out of any vertex, so there is
at least one vertex (let it be the vertex $j$) in the set $\W^F\setminus
\W^H$, which is not attainable in the forest $F$ from the set $\W^H
\setminus\W^F$, and hence the tree of the forest $F$, having the vertex
$j$ as a root, has not intersection with the set $\W^H\setminus\W^F$.
This way, all the vertices of the tree $T^F_j$ except the root $j$ itself,
belong to the set $(\N\setminus\W^F) \cap (\N\setminus\W^H)$, so arc goes
out of every vertex from the set $\V T^F_j\setminus\{ j\}$ in the forest
$H$ (and in $F$ naturally).

Let us construct preliminary forest $E\in\tfk $, which is necessary to the
final constructing of the descendant $G\in\tf^{k+1}$ of the forest $F$.
We take $\V T^F_j$-exchange of $F$ by $H$ as this auxiliary graph $E$, and
we designate $\V T^F_j$-exchange of $H$ by $F$ as $Q$. By force of
{Sequence 2} from  {Criterion} the graphs $E$ and $Q$ are forests.  The
forest $E$ contains one arc more than $F$, as there are no arc coming from
the vertex $j$ in $F$, but there is one in $H$ ($j\in\W^F\setminus\W^H$).
So $E\in\F^k$ and, analogically, $Q\in\F^{k+1}$ and

\be
\varphi_k \le\Sigma^E \ ,\ \ \varphi_{k+1}\le\Sigma^Q \ .  \la{1}
\ee
If we designate by $\Delta$ the quantity $\Sigma^E-\Sigma^F=\Sigma^E-
\varphi_{k+1}$, then, obviously,

\be
\Sigma^Q=\Sigma^H-\Delta \ .\la{2}
\ee
Using (\ref{1}) and (\ref{2}) we get $\Sigma^E=\varphi_{k+1}+\Delta\le
\Sigma^H -\Delta+\Delta=\varphi_k \ ,$ and hence $\Sigma^E=\varphi_k$,
what means that $E\in\tfk$.

Let the vertex $j$ in the forest $E$ belong to the tree $T^E_m$ with
vertex $m$ as a root. Consider the maximal walk being a subgraph  of the
tree $T^E_m$ and starting from the vertex $j$, all vertices of which
belong to the set of vertices of the tree $T^F_j$. Let $n$ be  final
vertex of this way.  Designate as $T$ maximal subtree of the tree $T^E_m$
with vertex $n$ as a root, all vertices of which belong to the set of
vertices of the tree $T^F_j$. Notice, that all trees of the forest $F$
with the exception of the tree $T^F_j$ are subtrees of the trees of the
forest $E$ with the same roots, but the vertices of the set  $\V T^F_j$
are "divided" among the trees of the forest $E$. So, we can confirm that
there are no arcs coming into the set $\V T$ in the forest $E$ and by
force of {Sequence 3} from {Criterion} graph $G$ being $\V T$-exchange of
$F$ by $E$ is a tree, and obviously it belongs to the set $\F^k$. If we
consider $\V T$-exchange of $E$ by $F$, which by force of the same
{Sequence 3} from {Criterion} is a thee, analogically to the previous we
are convinced that really $G\in\tfk$, but by the construction it is a
descendant of the forest $F$. To the other side the affirmation of the
theorem is proved analogically.

The theorem on "relatives" lets us easy prove known system of convexity
inequalities  \cite{VF,V}. Exactly, it is valid

{\bf Proposition 2.} {\it The quantities $\varphi_k$ satisfy to the
following chain of convexity in\-equali\-ties}

\be
\varphi_{k-1}-\varphi_k\ge\varphi_k- \varphi_{k+1} \ . \la{convex}
\ee

{\bf Proof.} By the theorem on "relatives" any forest $H\in\tf^{k-1}$ can
be constructed using redirection of arcs coming from the vertices of the
only tree of some forest $G\in\tfk$, which, in its turn, can be
constructed by redirection of arcs coming from the vertices of the only
tree of some forest $F\in\tf^{k+1}$.  Let $F$, $G$, $H$  be just such
"relative" forests. Then there is at least one tree of the forest $F$,
from every vertex of which an arc goes out in the forest $H$. Let the
vertex $i$ be root of this tree and let us designate by $f$ the sum of
weights of arcs coming in forest $F$ from the vertices of the set $\V
T^F_i$, and by $h$ -- the sum of weights of arcs outgoing from the
vertices of the same set in $H$.

Let $P$ be $\V T^F_i$-exchange of $F$ by $H$, and $Q$ be $\V T^
F_i$-exchange of $H$ by $F$. By force of the {Sequence 2} from the
{Criterion} both these  graphs are forests and belong to the set $\F^k$
(because there is not an arc coming from the vertex $i$ in the forest $F$,
but there is one coming from this vertex in forest $H$) and, hence

$$\Sigma^P=\Sigma^F+h-f=\varphi_{k+1}+h-f
\ge\varphi_k \ ,$$

$$\Sigma^Q=\Sigma^H-h+f=\varphi_{k-1}-h+f
\ge\varphi_k \ .$$
The Proposition is a direct sequence of the last two inequalities.

Note, that the following inequalities

\be
\varphi_{n-i}-\varphi_n\ge\varphi_{m+i}-\varphi_m \ , \ \
m\ge n \ , \ \min (N-m,n)\ge i\ge 0 \ , \la{convex2}
\ee
are the sequences from the system of convexity inequalities
(\ref{convex}).

Let us prove the following auxiliary

{\bf Proposition 3.} {\it Let $F\in\tf^n$ and $G\in\tf^m \ , \ \ m\ge
n$, and let $\D $ be subset of the set of vertices $\N $, such that graphs
$P$ and $Q$, being $\D$-exchange of $F$ by $G$ and $\D$-exchange of
$G$ by $F$ correspondingly, are forests.  Then if

a) $\D$ contains $l\ge 0$ roots of the forest $F$ more than roots of the
forest $G$, then $P\in\tf^{n-l}$ and $Q\in\tf^ {m+l}$;

b) $\D $ contains $l\ge m-n$ roots of the forest $G$ more than roots of
the forest $F$, then $P\in\tf^{n+l}$ and  $Q\in\tf^{m-l}$.}

{\bf Proof.} We prove point b) (point  ) can be proved analogically).
Designate as $\Delta$ the following quantity $\Delta=\Sigma^P-\Sigma^F=
\Sigma^G-\Sigma^Q \ .$ It is followed from the condition, that $P\in
\F^{n+l}$ and $Q\in\F^{m-l}$, so

$$\Sigma^P=\varphi_n+\Delta\ge\varphi_{n+l} \ ,\ \  \Sigma^Q=\varphi_m-
\Delta\ge\varphi_{m-l} \ .$$
Combining these two inequalities one gets $\varphi_{m}-\varphi_{m-l}\le
\varphi_{n+l}-\varphi_n \ . $ However from (\ref{convex2}) under $m\le
n+l$ it is followed reverse inequality and hence $\Sigma^P=\varphi_{n+l}$
and $\Sigma^Q= \varphi_{m-l}$ and this proves the proposition directly.

\section{Algebras of subsets}

At the present paragraph we will construct the system of embedded algebras
$\alk , \ k=1,2,\cdots N,$ of subsets of the set of all vertices $\N $ and
investigate the properties of the elementary sets of these algebras.

Let us consider all connected components $T$ (trees) of the forests
$F\in\tfk$. The sets of vertices of the trees $T$ are the base of the
algebra $\alk$ (i.e.  algebra $\alk$ is generated by the sets of vertices
$\V T$ of the trees of the forests $F\in\tfk$).

{\bf Theorem 2.} {\it The sequence of algebras $\alk$  is an
increasing one:

$$\{\N ,\emptyset\}=\al_1 \subseteq \al_2
\subseteq \cdots \subseteq  \al_{N-1} \subseteq \al_N =2^\N  \
,$$
where $2^\N $ is the set of all subsets of the set $\N $.}

{\bf Proof.} Direct sequence of {Theorem 1}.

Let us give a {definition}. We call the vertex $j$ as {\it marked point}
({\it vertex}) {\it of the level} $k$, if there exists at least one forest
$F\in\tfk$, where $j$ is a root (i.e. there exists connected component
$T^F_j$).

Elementary sets of algebras $\alk$ can as contain as not contain marked
vertices. Elementary set can contain few marked vertices at once. Those
elementary sets, that contain marked vertices we will call {\it marked
sets}.

Let $\S$ be some subset of the set of vertices $\N$. As $\tfk |_{\cal S}$
we will designate the set of subgraphs of the set of forests $\tfk$
induced by the set $\S$.

Let us see what the properties of extreme forests are in case, if under
some $k$ there is equality in the system of convexity inequalities
(\ref{convex}):

\be
 \varphi_{k-1}-\varphi_k =\varphi_k- \varphi_{k+1} \ . \la{equal}
\ee

{\bf Theorem 3.} {\it Let (\ref{equal}) be fulfilled, then

1) $\alk=\al_{k+1}$,

2) $\tf  ^{k-1}|_\E \subseteq \tfk |_\E
\supseteq\tf  ^{k+1}|_\E$, where $\E$ is an arbitrary elementary set of
the algebra $\alk$.}

{\bf Proof.} According to the {Theorem "on relatives"} every forest
$H\in\tf^{k-1}$ possesses at least one ancestor $F\in\tfk$, which in its
own, possesses at least one ancestor $G\in\tf^{k+1}$. Let $H$, $F$ and
$G$ be such relative forests. There are 2 possible scenarios of getting
granddescendant $H$ from grandancestor $G$. It is easy to see, that by one
of them 4 trees of the forest $G$ participate in the construction of the
forest $H$, and by another one -- only 3. Let us see on the first possible
scenario.

So, let $T^G_i$, $T^G_j$, $T^G_l$ and $T^G_m$ be trees of the forest $G$
with the roots $i$, $j$, $l$ and $m$ correspondingly. Let the forest $F$ be
constructed from the forest $G$ by uniting trees $T^G_i$ and $T^G_j$ with
may be redirecting of arcs coming from vertices of, for example, the
tree $T^G_j$, i.e. $T^F_i|_{\V T^G_i}=T^G_i \ , \ \ T^F_i|_{\V T^G_j}$ is
a tree and $\V T^F_i=\V T^G_j\cup\V T^G_i$, other trees of the forests $F$
and $G$ coincide between each other correspondingly. The forest $H$ in its
turn is received from the forest $F$ by uniting trees $T^F_l$ and $T^F_m$
with may be redirecting arcs outgoing from the vertices of, for example,
the tree $T^G_m$, i.e.  $T^H_l|_{\V T^F_i}=T^F_l \ , \ \ T^H_l|_{\V
T^F_m}$ is a tree and $\V T^H_l =\V T^F_m\cup\V T^F_l$, other trees of
the forests $F$ and $G$ coincide between each other correspondingly (note,
that also $T^G_l= T^F_l$ and $T^F_m=T^G_m$). Designate as $F' \ \ $ $\V
T^G_j$-exchange of $H$ by $G$. It is obvious (by {Sequence 2} from
{Criterion} and {Proposition 3}), that $F'\in\tilde \F^k$. By this every
tree of the forest $G$ and every tree of the forest $H$ is either a tree
of the forest $F$ or a tree of the forest $F'$, that confirms both points
of the theorem. Another variant of the scenario is considered
analogically.

We say, that the vertex $j$ is attainable from the vertex $i$ at the
level $k$ or simply $j$ is  $k$-attainable from $i$, if there is at least
one forest $F\in\tfk $, such as there is  $i\cdot j$-walk in
$F$.

Let us see what are the properties of extreme forests in case if under
some $k$ there is strong inequality in the system of convexity
inequalities:

\be
 \varphi_{k-1}-\varphi_k >\varphi_k- \varphi_{k+1} \ . \la{nonequal}
\ee

{\bf Proposition 4.} {\it Let (\ref{nonequal}) be taken place and $i$ and
$j$ be level $k$ marked vertices. Let also the vertex $i$ be attainable
from the vertex $j$ on the level $k$, then the vertex $j$ is attainable
from the vertex $i$ on the level $k$ and, moreover, the vertices $j$ and
$i$ belong to the same marked set of this level.}

{\bf Proof.} Under condition there is such forest $F\in\tfk $,
where the vertex $i$ is rear comparative to the vertex $j$. Without loss
of generality one can consider that, the vertex $j$ is a root in the
forest $F$ (otherwise, if some marked vertex $m$ is a root of the tree
containing the vertices $i$ and $j$ at this forest, the following
discussions one can lead for any pair of vertices $i$ and $m$ or $j$ and
$m$). Suppose, that there is such forest $G\in\tfk $, in
which the vertex $j$ is a root, and the vertex $i$ does not belong to the
tree having $j$ as a root. Let $\D  =\V T^F_i \cap \V T^G_j$, and $P$ and
$Q$ are  $\D$-exchanges of $F$ by $G$ and of $G$ by $F$ correspondingly.
Then by {Proposition 3} $ \ \ P\in\tf^{k+1}$ and $Q\in\tf^{k-1}$. Let us
denote by $f$ and $g$ the sums of weights of arcs coming from the vertices
of the set $\D $ at forests $F$ and $G$ correspondingly, then

$$\varphi_{k+1}=\Sigma^P=\Sigma^F-f+g=\varphi_k-f+g \ , $$

$$\varphi_{k-1}=\Sigma^Q=\Sigma^G+f-g=\varphi_k+f-g \ , $$
whence it follows that $\varphi_{k-1}-\varphi_k =\varphi_k- \varphi_{k+1}
\ , $ which contradicts (\ref{nonequal}). So, in any forest $G\in\tfk$, in
which the vertex $j$ is a root, the vertex $i$ belongs to the set of
vertices of the tree $T^G_j$. From here it easy follows, that there is not
such a forest in the set $\tfk$, in which the vertices $i$ and $j$ belong
to different trees, which means validity of the proving proposition.

Note, that this proposition means in particular that if (\ref{nonequal})
is fulfilled, so every marked set of algebra $\alk$ contains exactly
one root of an arbitrary forest $F\in\tfk$, and it is valid the following.

{\bf Theorem 4.} {\it Let (\ref{nonequal}) be fulfilled, then the algebra
$\alk$ contains exactly $k$ marked elementary sets.}

{\bf Proof.} Any forest $F\in\tfk$ consists of $k$ trees and hence, there
are not less than $k$ marked elementary sets in $\alk$.  These $k$ marked
sets are those elementary sets that contain the roots of the trees of $F$.
Any root of an arbitrary forest $G\in\tfk$ naturally belongs to one of
the trees of the forest $F$ and, hence, some root of the forest $F$ is
accessible from it (root of $G$), and it means by Proposition 4 that this
root belongs to one of mentioned elementary sets. Thus, there are exactly
$k$ marked sets in $\alk$.

Let us call as {$k$-attraction domain} of marked vertex $i$ such set of
vertices, which consists of such vertices $j$  that $i$ is accessible from
$j$ in at least one forest $\F\in\tfk$.

{\bf Proposition 5.} {\it Let (\ref{nonequal}) be fulfilled for some $k$,
then for every marked vertex $i$ there is such forest $\F\in\tfk$, in
which the vertex $i$ is a root and the set of vertices of the tree $T^F_i$
coincides with ${k}$-attraction domain of the vertex $i$, and also the
sets of ${k}$-attraction domains of mutually ${k}$-attainable vertices
coincide with each other.}

{\bf Proof.} Let $F$ and $G$ be forests belonging to the set $\tfk$, in
which mutually $k$-attainable vertices $i$ and $j$ (in particular they can
coincide) are roots of the trees $T^F_i$ and $T^G_j$ correspondingly.  It
is sufficient to show, that there is such a forest $H\in\tfk$, where
the vertex $i$ is a root and $\V T^H_i \supseteq \V T^F_i \cup \V T^G_j$.
Let $\D = \V T^G_j\setminus\V T^F_i$, then by Proposition 3  $\D$-exchange
$F$ by $G$ is required forest $H$.

{\bf Proposition 6.} {\it Let (\ref{nonequal}) be fulfilled, $F\in\tfk$
and $\E$ is elementary set belonging to algebra $\alk$, then there is
such forest $G\in\tfk$, where all arcs coming out from the vertices
of the set $\E$, coincide with ones coming out from them in the forest
$F$, and also there are no arcs coming into the set $\E $ from the outside
in $G$.}

{\bf Proof.} Let there be an arc coming into the set $\E $ from some
elementary set $\E _1$ in the forest $F$. Since the sets $\E$ and $\E _1$
are elementary, so there is such forest $H\in\tfk$, where both these sets
belong to different trees.  Let $\E$ belong to the tree with $i$ as a root
in $F$, and $\E _1$ belong to the tree with $j$ as a root in the forest
$H$.  Let $\D$ be the set $\V T^F_i\cap\V T^H_j$. Let $G$ be $\D$-exchange
$F$ by $H$. By Proposition 4 the vertices $i$ and $j$ simultaneously
belong or do not belong to the set $\D$. So $G\in\tfk$ and there are not
any arcs coming into the set $\E$ from the set $\E _1$ in this forest, and
also there are not more additional arcs coming into the set $\E$ in $G$,
in comparison to ones coming into $\E$ in the forest $F$. If there are
some arcs coming into the set $\E$ in $G$, one can repeat the procedure
above now concerning the forest $G$ and get the forest, where no one arc
comes into the set $\E$, but all arcs coming from it coincides with those
coming from vertices of $\E$ in the forest $F$.

Next proposition being direct consequence of Proposition 6 is in some
sense inverse to Proposition 5. If Proposition 5 tells how big tree of
extreme forest can be, but in the following one we explain how small it
can be.

{\bf Proposition 7.} {\it Let (\ref{nonequal}) be fulfilled for some $k$,
then for every marked elementary set $\M$ of algebra $\alk$ there is such
forest $F\in\tfk$, where $\M$ is a set of vertices of one of trees of
$F$, and also there is not such a forest belonging to $\tfk $, where arcs
come out of the set $\M$. }

{\bf Proof.} Let us suppose inverse. Let $F\in\tfk$ be a forest, where at
least one arc comes out of $\M$ with, let us say, the vertex $m$ as an
origin. By Proposition 6 without loss of generality one can suppose that
there are not any arcs coming into $\M$ from outside. In addition,
according to Proposition 4, the set $\M$ contains exactly one root of
$F$. But then the tree of $F$ having this root does not contain the vertex
$m$ and is contained in $\M$, which is in contradiction with elementary
character of $\M$.

Proposition 7 means in particular, that any subgraph of an arbitrary
forest $F\in\tfk$, induced by marked elementary set of algebra $\alk$, is
a tree if (\ref{nonequal}) is fulfilled.  It is prove to be that indicated
property is valid for unmarked elementary sets too.

{\bf Theorem 5.} {\it Let (\ref{nonequal}) be fulfilled, then induced by
any elementary set $\E$ of algebra $\alk$ subgraph of any forest
$F\in\tfk$ is a tree.}

{\bf Proof.} It is necessary to show, that not more than one arc can come
out of an arbitrary elementary set $\E $. Let $F\in\tfk$. According to
Proposition 6 one can suppose, that there are not any arcs coming into
$\E$ from outside in $F$. Let us verify firstly, that not more than one
arc can come out of the set $\E$ into any other elementary set. On the
contrary, we assume that there are, for example, two arcs at the forest
$F\in\tfk$ coming out of the set $\E$ into some elementary set $\E _1$
of algebra $\alk$. Let also the arcs coming out of the set $\E$ into
$\E_1$ have their origin at the vertices $a$ and $b$ and let the sets $\A$
and $\B$ be sets of rear vertices with respect to vertices $a$ and $b$
correspondingly (including vertices $a$ and $b$ themselves). The sets $\A$
and $\B$ do not intersect with each other and $\A\cup\B =\E$. As the
sets $\E$ and $\E _1$ are elementary, so there is such forest $G\in\tfk$,
where these sets belong to different trees, let us say, to the trees
$T^G_j$ and $T^G_m$ correspondingly. Let $H$ be $\A$-exchange of $G$ by
$F$. Obviously, that $H\in\tfk$. In addition, since $\E$ is elementary
and, hence, its vertices at any forest from the set $\tfk$ must belong to
the same tree, among them at $H$ too. It is possible only if the vertices
of the set $\B $ are rear with respect to the vertex $a$ at the forest $G$
(only in this case elementary set $\E $ belongs entirely to single tree at
$H$, namely to the tree $T^H_m$).  Analogously, if $Q$ is $\B$-exchange
of $G$ by $F$, so $Q\in\tfk$ and the vertices of the set $\A$ must be
rear with respect to the vertex $b$ at the forest $G$. So the vertices $a$
and $b$ are rear with respect to each other at $G$, which is impossible
because $G$ is a forest.

Other cases, where arcs could come out of $\E$ into several elementary
sets one can examine analogously.

{\bf Theorem 6.} {\it Let (\ref{nonequal}) be fulfilled, then

i) induced by any elementary set $\E$ belonging to the algebra $\alk$
subgraph of an arbitrary forest $F\in\tf^{k-1}$ is a tree,

ii) if $\U$ is unmarked elementary set belonging to the algebra $\alk$,
then $\tf^{k-1}|_\U  = \tfk |_\U$.}

{\bf Proof.} Let $F$ and $G$ be relative forests belonging correspondingly
to $\tfk$ and $\tf^{k-1}$, and let also one can construct the forest $G$
from $F$ by adding an arc coming out of the root $i$ of some tree $T^F_i$,
and, may be, by redirecting of arcs that come out of other vertices of
this tree. According to Proposition 6, without loss of generality, one can
consider that $\M =\V T^F_i$ is marked elementary set, and by theorem on
"relatives" the graph $G|_\M$ is a tree. In this case, if $\U$ is
unmarked elementary set of the algebra $\alk$, so $G|_\U =F|_\U $.

Theorems 5 and 6 are very important for the consequent constructions,
since based on Theorem 5 one can construct enlarged graphs and to
determine adjacencies (and their weights) connecting enlarged vertices
(elements of decomposition of the set of all vertices). Theorem 6
allows based on one level of enlargement to construct the following one.

\section{Weighted condensations}

Proved above properties of extreme forests allow us to look on them and
at all on directed graphs with weighted adjacencies in "an enlarged way",
without interest on details of their arc connections inside elementary
sets, but paying attention only on connections among elementary sets,
understanding elementary sets themselves as a vertices of some enlarged
graph. Let us convert what has been said above into precise definition.
Beforehand we remind existing definition of condensation for non-weighted
directed graph, which just allows understand graphs in an enlarged way.
Here is the corresponding definition.

Let $\{ {\cal S}_1,{\cal S}_2,\cdots ,{\cal S}_M\} $ be strong components
(strong component is the set of inter-attainable vertices) of digraph $G$.
{\it Condensation} of digraph $G$ is  digraph $\hat G$ with the set of
vertices $\{ s_1,s_2,\cdots ,s_M\} $, where the pair $(s_i,s_j)$ is an arc
in $\hat G$ if and only if there is an arc in $G$ with origin belonging to
${\cal S}_i$, and terminus belonging to ${\cal S}_j$.

Mentioned definition is rather poor, since, for example, for strong
digraphs (where all vertices are inter-attainable) condensation is trivial
and consists of only one vertex, and hence, there are not any arcs in it.
So we essentially modify the concept of condensation for weighted digraphs.

Let us firstly consider the case of non-directed graphs. In some sense the
following simple theorem \cite{E} is more strong reformulation of Theorem
on relatives but for non-directed graphs.

{\bf Theorem 7.} {\it Let the edge $e$ of non-directed graph $P$ possesses
the minimal weight among all edges, in which exactly one endpoint belongs
to the tree $T$ which is subgraph of $P$. Then there is at least one
spanning tree containing $T\cup e$ and having minimal weight  among
all spanning trees of $P$ containing $T$.}

According to this theorem all examinations drawn are valid but essentially
simplify. For example the division on marked and unmarked sets vanishes
(every set is marked) and also there is no necessity  to replace edges
under joining trees as it was in case of directed graphs (one only need
add an edge to connect two trees).  Of course for non-directed graph $P$
with weighted adjacencies inequalities of convexity are fulfilled and if
(\ref{nonequal}) is valid then algebra $\alk$ contains exactly $k$
elementary sets. The main property resulting from this theorem, that is
useful for us, we point out as following.

{\bf Property 1.} {\it Subgraph of any forest $F\in\tf^n$ induced by
elementary set $\E$ of algebra $\alk$, $n\le k$, is a tree.}

{\bf Property 2.} {\it For every forest $F\in\tf^n$ there exists such
forest $G\in\tfk$, $n\le k$ (and, into opposite side, for any $G\in\tfk$
there exists such $F\in\tf^n$ ) that $F|_\E =G|_\E$, where $\E$ is an
arbitrary elementary set of algebra $alk$.}

{\bf Definition.}  Let $P$ be non-directed graph with weighted
adjacencies $p_{ij}$, and let (\ref{nonequal}) be fulfilled, $\alk$
-- algebra of subsets of the set of all vertices generated by the sets of
vertices of trees belonging to $\tfk(P)$. We call non-directed graph $P^k$
with $k$ vertices as {\it weighted condensation of the level} $k$ (simply
-- {\it k-weighted condensation)} of $P$ if weights of it adjacencies are
equal to the following numbers

\be
p^k_{xy}=\min\limits_{\sc i\in\X  \atop
\sc j\in\Y } p_{ij}\ , \la{vadn}
\ee
where $\X$ and $\Y$ are elementary sets of algebra $\alk$. If there is not
any edge in $P$, such that one of its ends belongs to the elementary set
$\X$, and another to the elementary set $\Y$, so we suppose that there is
not corresponding edge $(x,y)$  in $P^k$.

It seems natural to consider that in graph of weighted condensation not
only arcs possess weights but vertices too, which are actually elementary
sets of corresponding algebra. We determine weight of vertex $s$, or which
is the same, weight of elementary set $S$ corresponding to vertex $s$, as
minimum of weight of spanning tree of graph $P|_{\cal S}$, i.e. as the
quantity

$$\mathop{\min}\limits_{\sc T\subset P \atop
\sc \V T={\cal S}} \sum\limits_{(i,j)\in T}p_{ij} \ .
$$
As weighted condensations, represent themselves usual graphs with weighted
adjacencies, so all previous properties are valid for them (introduction
of weights of vertices is not change anymore because we consider only
spanning subgraphs, which include all vertices by definition).  In
particular, one can consider factor-forests of $P^k$ and determine the
sets $\F^n(P^k)$ and also their subsets $\tf^m(P^k)$ possessing minimal
weight. The weight itself of the forest $F\in\tf^n(P^k)$ we determine
as stated above in the following way

\be
\Sigma^F=\sum\limits_{(x,y)\in F}p^k_{xy}+ \sum\limits_{\E \in
\alk} \mathop{\min}\limits_{\sc T\subset V \atop \sc
\V T=\E } \sum\limits_{(i,j)\in T}p_{ij} \ , \la{SFN}
\ee
where $\E $ is elementary set of algebra $\alk$. For example, any
forest $F\in\F^k(P^k)$ is empty graph ($k$ vertices (however possessing
their own weights) and no edges), any $F$ belonging to $\F^1(P^k)$  is a
spanning tree of $P^k$.

Under definition (\ref{SFN}) it is obvious that if we introduce the
numbers $\varphi^k_n$, $n\le k$, by the rule

\be
\varphi^k_n =\mathop{\min}\limits_{F\in\F ^n(P^k)}\Sigma^F \ ,
\la{defvarphikn}
\ee
then by force of Property 1
\be  \varphi_n = \varphi_n^k \ , \ \  n\le k \ ,
\la{varphi=}
\ee
and, of course, inequalities of convexity are valid:

\be
\varphi_{n-1}^k-\varphi_n^k\ge \varphi_n^k-\varphi_{n+1}^k \ , \ \
n=2,3,\cdots ,k-1 \ .  \la{nonequalkn}
\ee
Equalities (\ref{varphi=}) mean exactly, that minimum weights of spanning
trees, consisting of equal number of trees $n\le k$, of weighted
condensation $P^k$ and graph $P$ proper coincide with each other.

Now we consider analogical examination for directed graphs. Let $V$ be
digraph with weighted adjacencies $v_{ij}$, and let (\ref{nonequal}) be
fulfilled, $\alk$ -- algebra of subsets of the set of vertices of
$V$, generated by the sets of vertices of trees of forests belonging to
$\tfk(V)$. Algebra $\alk$ contains at least $k$ elementary sets, to be
precisely, it contains $k+l$ elementary sets, where $l$ is the number of
unmarked sets (this number can be equal to zero).

{\bf Definition.} Let us call digraph $V^k$ with $k+l$ vertices as {\it
weighted condensation of the level} $k$ (simply -- {\it k-weighted
condensation)} if weights of its adjacencies are equal to the following
numbers

\be
v^k_{xy}=\min\limits_{\sc i\in\X  \atop
\sc j\in\Y } (\min\limits_{T_i(\X )}\Sigma^{T_i({\cal
X})}+v_{ij})\ , \la{vad}
\ee
where $\X $ and $\Y $ are elementary sets of algebra $\alk$, $T_i(\E)$ is
a tree with $\E$ as a set of vertices and $i$ as a root. If there is not
any arc in $V$, such as its origin belongs to the elementary set $\X$,
and the terminus  to the elementary set $\Y$, and under this  $i$ is a
root of at least one spanning tree of digraph $V|_\X$ we suppose that
there is not arc $(x,y)$ in $V^k$.

The necessity of weights determination in a different way than it was in
non-directed situation is caused by the fact that one must be sure that
the set $\Y$  is attainable from every vertex of $\X$  and in this case
only it is justified to introduce an arc $(x,y)$ into  graph $V^k$.  Note,
that weight minimum of tree $T_i(\X)$ depends on vertex $i$, so generally
speaking in the situation of directed graphs it is not possible to
introduce the weight of elementary set and one needs  add "it" (look at
(\ref{vad})) to corresponding arc going out of this set.  Nevertheless, if
there are not arcs going out of some set $\X$ in digraph, it is possible
to determine weight of $\X$ as minimum by all $i\in\X$ of weights of
trees $T_i(\X)$.

As graph $V^k$ has at least $k$ ($k+l$ to be precisely) vertices one can
consider, in particular, spanning forests of it and to determine the sets
$\F^m(V^k)$ and also their subsets $ \tf^m(V^k)$ possessing minimal
weight. However the weight itself of the forest $F\in\tf^m(V^k)$ we must
determine in other way than in non-directed situation, because arc weights
(\ref{vad}) are determined not analogous to edge ones (\ref{vadn}).
Namely:

\be
\Sigma^F=\sum\limits_{(x,y)\in F}v^k_{xy}+ \sum\limits_{\sc \E\in
\alk \atop \sc d^+(e)=0} \mathop{\min}\limits_{ \sc  T\subset V \atop \sc
\V T= \E } \sum\limits_{(i,j)\in T}v_{ij} \ , \la{SF}
\ee
where $F\in\F^m(V^k)$, $e$ is a root of $F$ corresponding to the
elementary set $\E\in\alk$. So weight of $F\in\F ^m(V^k)$ is determined
as sum of all arc weights $v^k_{xy}$ plus "weights" of those elementary
sets of algebra $\alk$, corresponding to which vertices in $F$ are roots.

Now one can introduce the quantities $\varphi^k_n \ , \ n=1,2,\cdots ,k$
by the rule analogous to (\ref{defvarphikn})

$$
\varphi_n^k=\min\limits_{F\in\tf^n(V^k)}\Sigma^F \ .
$$
and, of course, for these quantities the inequalities of convexity
(\ref{nonequalkn}) continue to be fulfilled, but (\ref{varphi=}) is not
true now and one can assert only that

$$\varphi_n\le \varphi_n^k \ , $$
as the minima $\varphi_m$ are calculated using graph $V$ itself, but the
numbers $\varphi_m^k$ -- only using its weighted condensation. However, by
force of definition of weighted condensations and its adjacencies
(\ref{vad}) $\varphi_k^k=\varphi_k$. Moreover, since by Theorem 6 subgraph
of any graph belonging to $\tf  ^{k-1}$ induced by an arbitrary
elementary set of algebra $\alk$ is a tree, so $\varphi_{k-1}^k=
\varphi_{k-1}$. Note, that (\ref{varphi=}) is a sequence of Property 1,
which is not valid here generally speaking.

Point here  that one can use the definition of weighted condensations
in case of non-fulfillment of (\ref{nonequal}) also. Namely, let under
some $k$ and $n\le k-1$

\be
\varphi_{k-n-1}-\varphi_{k-n}>\varphi_{k-n}-\varphi_{k-n+1}=\cdots =
\varphi_{k-1}-\varphi_k >\varphi_k-\varphi_{k+1} \ , \la{kn}
\ee
then, as it follows from Theorem 3, algebras $\al_{k-n+1}$, $\al_{
k-n+2},\ldots ,\alk$  coincide with each other and so the definition
of weighted condensations, initially introduced for index equal to $k$,
one can spread to indices $k-1$, $k-2$, $\ldots $, $k-n+1$. Under that it
is obvious that all this condensations are the same, so the number of
different condensations equal to the number of sign $'>'$ at the system of
convexity inequalities (\ref{nonequal}) plus one. Theorems 3 and 6 mean
also that under (\ref{kn}) subgraph of any forest belonging to one of sets
$\tf^{k-l}$, $l=1,2,\ldots ,n$, induced by arbitrary elementary set of
algebra $\alk$, is a forest and hence

$$\varphi_{k-l}=\varphi_{k-l}^k \ , \ \ l=1,2,\cdots ,n \ .$$
One could think that (\ref{varphi=}) is valid for digraphs, however it is
not so, because under (\ref{kn}) one has not any reason to expect that
subgraph of $F\in\tf^{k-n-1}$, induced by elementary set of $\alk$, is a
forest (and really it is not so, one can easy construct such example).
Nevertheless (\ref{varphi=}) takes place if the adjacencies $v_{ij}$ of
digraph $V$ can be written in the form

\be
v_{ij}=p_{ij}-p_{ii} \ , \la{pot}
\ee
where the numbers $p_{ij}\in R^1$ are weights of edges of some
non-directed graph $P$ ($p_{ij}=p_{ji}$). This property we will call as
{\it potentiality } of weights of digraph $V$. Such definition is bound up
with the fact, that under fulfillment of (\ref{pot}) the weights $v_{ij}$
can be realized as potential bars necessary to overtake in order to get
into point $j$ from point$i$ (the number $p_{ij}$ is transition potential
from $i$ to $j$, $p_{ii}$ -- potential of point $i$).  Equalities
(\ref{varphi=}) succeed from the following

{\bf Theorem 8.} {\it Let digraph $V$ possess potential weights
and its adjacencies satisfy (\ref{pot}), then Property 1 is valid for $V$
and (\ref{varphi=}) takes place.}

{\bf Proof.} From the definition of potentiality it is followed that if
there is an arc $(i,j)$ in digraph $V$, so there is an opposite arc
$(j,i)$ there. Further, for potential graph it is not difficult to see
that if some $i\cdot j$-way possesses minimum weight (minimum sum of arc
weights (potential bars)) among all ways from $i$ to $j$, then if one
changes all these arcs to opposite ones in this $i\cdot j$-way, one gets
$j\cdot i$-way with minimum weight among all ways from $j$ to $i$ in $V$.
Now let us turn to Theorem 1 (on "relatives"). According to it any forest
$G\in\tf^{k+1}$ one can construct from some forest $F\in\tfk$ by adding an
arc connecting two trees and may be by redirecting of arcs in that tree,
from which this additional arc would go out. Let $F$ and $G$ be such
relative forests, and let $G$ one can get from $F$ by adding arc $(i,m)$,
where $i$ belongs to the set of vertices of tree $T^F_j$ with $j$ as a
root, and ,it is clear, if $i$ does not coincide with $j$, by
reconfiguration of arcs of this tree in such a way as to get on the set
$\V T^F_j$ a new tree, but with $i$ as a root. In this connection this
tree $G|_{\V T^F_j}$ must possess minimum weight among all trees on the
set $\V T^F_j$ with $i$ as a root. Let us construct new tree $G'\in\tfk$
from $F$ by adding the same arc $(i,m)$, but reconfiguration of arcs of
$T^F_j$ will be done in the following manner. Consider $i\cdot j$-way
belonging to tree $T^F_j$. It is, of course, the only in this tree and it
possesses minimum weight among all $i\cdot j$-ways in induced subgraph
$V|_{\V T^F_j}$.  Now change in $F$ arcs of this $i\cdot j$-way into
opposite ones (one gets under this a tree on the set $\V T^F_j$ with
minimum weight among all trees on this set with $i$ as a root) and add arc
$(i,m)$. This forest let call $G'$. It is extreme, of course, because it
was constructed under really minimum changes of forest $F$. Note, that
this forest $G'$ possesses one important property. If one takes away the
orientation from $F$ and $G'$, then these graphs coincide with each other,
except adjacency $(i,m)$ proper. It appears from the above the validity of
Property 1 for potential digraphs.

Theorem 8 shows, that the analysis of potential digraphs is not more
difficult than the same of non-directed graphs, and for them instead of
Property 2 it is valid

{\bf Property 2'.} {\it Let $V$ be potential digraph, then  for any
forest $F\in\tf^n(V)$ there exist such forest $G \in \tfk(V)$, $n\le k$
(and for any $G \in \tfk (V)$ there is such $F \in\tf^n(V)$), that induced
by any elementary set $\E \in\alk$ subgraphs of $F$ and $G$ coincide with
each other to within the orientation.}

So, our considerations above mean the following. Let us suppose that we
constructed $k$-weighted condensation of some directed graph $V$ and try
to build up condensation $V^n$, $n<k$, of some next level $n$. The
question appears: Can one do it using the information on already
constructed condensation only? It turns out that can not, generally
speaking.  More exactly, one can construct the corresponding algebra
$\al_n$, but  new adjacencies -- can not. One have to use information
on arcs of the initial digraph $V$. Nevertheless, if weights of $V$ are
potential  (or, moreover, graph is non-directed), it is not necessary to
use any additional information and to realize the transition to next
hierarchy level one can forget "prehistory" of graph and use the
adjacencies of $V^k$ only. This reduces considerably the number of
calculations required.

\section{Instead of discussion}

The method suggested can have a lot of applications in different brunches
of science such as economy and finances, biology and neuron-nets,
probability theory and random processes, mathematical and theoretical
phy\-sics.  This is forced just by the necessity to determine the
structure and the hierarchy of complicated objects and using this
information to give a conclusion which processes are essential on each
level and which are not. For example, at exponentially large times in
dynamic systems under small random perturbations some sublimit
distributions appear \cite{VF}. They correspond in fact to distributions
concentrated at marked elementary sets of some algebra $\alk$, the
number of nontrivial possible time scales is equal to the number of
different algebras. Under this, the generators of Fokker-Plank type
equations (being singular perturbed ones \cite{L}), which govern
distribution functions of stochastic differential equations, possess
very special spectrum. Its low-frequency spectrum and corresponding
eigenfunctions are determined by weighted condensations of some
special digraph \cite{V1,V2}, which analysis connects with the opportunity
of representing of characteristic polynomial in terms of tree-like structure
of corresponding digraph \cite{CDZ,FS,V3}.

{\rm This work was  supported RFBR, grants N-99-01-00696 and
N-98-01-01063.} \vskip.5cm


\end{document}